\documentclass[11pt,dvips,twoside]{article}
\usepackage[usenames]{color}

\usepackage{pslatex}
\usepackage{fancyhdr}
\usepackage{graphicx}
\usepackage{geometry}
\usepackage{times}
\usepackage{mathptmx}
\usepackage{array}
\usepackage{moreverb}
\usepackage{url}

\def\figurename{Figure} 
\makeatletter
\renewcommand{\fnum@figure}[1]{\figurename~\thefigure.}
\makeatother

\def\tablename{Table} 
\makeatletter
\renewcommand{\fnum@table}[1]{\tablename~\thetable.}
\makeatother

\usepackage{amsmath}
\usepackage{amssymb}
\usepackage{amsfonts}
\usepackage{amsthm,amscd}

\newtheorem{theorem}{Theorem}[section]

\newtheorem{proposition}[theorem]{Proposition}
\theoremstyle{definition}

\theoremstyle{remark}

\numberwithin{equation}{section}

\def\R{\mathbb R}

\begin{document}

\title{\bfseries\scshape{Asymptotic stability of solutions to abstract 
differential equations}}
\author{\bfseries\scshape A. G. Ramm\thanks{e-mail address: \tt{ramm@math.ksu.edu}}\\
Kansas State University, Department of Mathematics,\\
Manhattan, KS 66506-2602, USA
\\ \\ \\{\rm (Communicated by Editor-in-Chief)}}

\date{}
\maketitle

\thispagestyle{empty} \setcounter{page}{1}

\thispagestyle{fancy} \fancyhead{}

\fancyhead[L]{\large {\LARGE J}ournal of {\LARGE A}bstract {\LARGE D}ifferential {\LARGE E}quations and {\LARGE A}pplications\\
Volume 1, Number 1, pp. {\thepage--\pageref{lastpage-01} (2010)}}
\fancyhead[R]
\fancyfoot{}
\renewcommand{\headrulewidth}{.0000pt}

\noindent\hrulefill

\noindent {\bf Abstract.} An evolution problem for abstract
differential equations is studied. The typical problem is:
$$\dot{u}=A(t)u+F(t,u), \quad t\geq 0; \,\, u(0)=u_0;\quad \dot{u}=\frac
{du}{dt}\qquad   (*)$$ Here $A(t)$ is a linear bounded operator in a
Hilbert space $H$, and $F$ is a nonlinear operator, $\|F(t,u)\|\leq
c_0\|u\|^p,\,\,p>1$, $c_0, p=const>0$. It is assumed that
Re$(A(t)u,u)\leq -\gamma(t)\|u\|^2$ $\forall u\in H$, where
$\gamma(t)>0$, and the case when $\lim_{t\to \infty}\gamma(t)=0$ is
also considered. An estimate of the rate of decay of solutions to
problem (*) is given. The derivation of this estimate uses a
nonlinear differential inequality.

\noindent \hrulefill

\vspace{.3in}

\noindent {\bf AMS Subject Classification:} 34G20, 37L05, 44J05,
47J35.

\vspace{.08in} \noindent \textbf{Keywords}: Nonlinear inequality;
Asymptotic stability; Abstract differential equations.

\section{Introduction}
Let
\begin{equation}
\label{e1}
\dot{u}(t)=A(t)u+F(t,u), \quad t\geq 0; \quad \dot{u}(t):=
\dot{u}=\frac{du}{dt},
\end{equation}
\begin{equation}
\label{e2}
u(0)=u_0,
\end{equation}
where $u\in H$, $H$ is a Hilbert space, $A(t)$ is a bounded linear
operator in $H$, $F(t,u)$ is a nonlinear operator,
\begin{equation}
\label{e3}
\|F(t,u)\|\leq c_0\|u\|^{p}, \quad p>1,
\end{equation}
$c_0$ and $p$ are positive constants, and $u_0\in H$.

One says that $A(t)\in B(\rho, N)$ if every solution to the equation
\begin{equation}
\label{e4}
\dot{v}(t)=A(t)v
\end{equation}
satisfies the estimate
\begin{equation}
\label{e5}
\|v(t)\|\le N e^{\rho (t-s)}\|v(s)\|, \qquad t\geq s\geq 0,
\end{equation}
where $N>0$ and $\rho$ are real numbers. This definition is
discussed in \cite{DK} and goes back to P. Bohl (see the historical
remarks in \cite{DK}). If $U(t,s)$ is an operator that solves the
problem
\begin{equation}
\label{e6}
\dot{U}(t,s)=A(t)U(t,s), \qquad t\ge s; \quad U(s,s)=I,
\end{equation}
 where $I$ is the identity operator, then \eqref{e5} is equivalent (by the
Banach-Steinhaus theorem) to the estimate
\begin{equation}
\label{e7}
\|U(t,s)\|\le N e^{\rho (t-s)}, \qquad t\geq s\geq 0.
\end{equation}
Let us define, following \cite{DK}, the notion of upper general
exponent $\kappa$ for the solutions to \eqref{e4}:
\begin{equation}
\label{e8}
\kappa={\overline{\lim}}_{t,s \to \infty}\frac{\ln \|U(t+s,s)\|}{t},
\qquad
t,s\geq 0.
\end{equation}

If $\kappa<0$, then $\|v(t)\|=O(e^{-|\kappa|t})$ as $t\to \infty$,
$s$ being fixed.

The following result is obtained in \cite{DK}, Theorem 3.1, Chapter
7.
\begin{proposition}
 {\it If $\kappa<0$ and assumption \eqref{e3}
holds, then the zero solution to equation \eqref{e1} is
asymptotically stable. }
\end{proposition}
Recall that the zero solution to equation \eqref{e1} is called
Lyapunov stable if for every $\epsilon>0$, one can find a
$\delta=\delta(\epsilon)>0$, such that if $\|u_0\|\le \delta$, then
the solution to problem \eqref{e1}$-$\eqref{e2} satisfies the
estimate $\sup_{t\ge 0}\|u(t)\|\le \epsilon$. If, in addition,
$\lim_{t\to \infty}\|u(t)\|=0$, then the zero solution to equation
\eqref{e1} is called asymptotically stable in the Lyapunov sense.

As one can see from our proof of Theorem 1.2, the condition of
smallness of the initial data $\|u_0\|\le \delta$ can be replaced by
a different condition:  if  $\|u_0\|$ is arbitrary fixed, then one
still derives the relation  $\lim_{t\to \infty}\|u(t)\|=0$ from
\eqref{e23} (see below), provided that $c_0$ is sufficiently small.

In Proposition 1.1, the exponent $\kappa<0$ is a constant. For
example, if $A(t)=A^*(t)$ is a selfadjoint compact operator, and
$\lambda_j(t)$ are its eigenvalues, $\lambda_j(t)\leq
\lambda_m(t)<0$ if $j>m$, $j=1,2,3,...$, then $\lambda_1(t)\le
\kappa<0$.

Our goal is to derive an analog of Proposition 1.1
 such that $\lim_{t\to \infty}\lambda_1(t)=0$ is allowed, that is,
we do not assume that the spectrum $\sigma(A(t))$ of $A(t)$ lies in
a half-plane Re $z\le \kappa$, where $\kappa<0$ is a {\it fixed}
constant independent of $t$.

It is known (see, e.g.,  \cite{DK}) that if $A$ is a bounded linear
operator in $H$ with the spectrum $\sigma(A)$, which lies in the
half-plane Re $z\le -|\kappa|$, $|\kappa|>0$, then there is a
positive-definite operator $W$ such that Re$WA=-V$, where $V$ is an
arbitrary given positive-definite operator in $H$. In other words,
if $m=const>0$, $\sigma(A)\subset \{z: Re z\leq -|\kappa|<0\}$ and
$V=V^*\geq m>0$,  that is, $(Vu,u)\ge m(u,u)$ $\forall u\in H$, then
the operator equation $A^*W+WA=-2V$ is solvable for $W$. In fact,
there is an explicit formula for $W$: $W=2\int_0^\infty
e^{A^*t}Ve^{At}dt$  (see \cite{DK}). By Re$A$ one understands the
operator defined by the formula Re$A:=A_R:=(A+A^*)/2$, and
$A=A_R+iA_I$, where $A_R$ and $A_I$ are selfadjoint operators that
are called  the real and  imaginary parts of $A$. If $A_R\le -a$,
$a=const>0$, then $\sigma(A)$ lies in the half-plane Re $z\le -a$.
The notation $A\le -a$ means $(Au,u)\le -a(u,u)$ $\forall u\in H$.

The converse is not true: it is not true that if the spectrum of a
linear bounded operator $A$ lies in the half-plane Re $z\le -a$,
then the inequality $A_R\le -a$ holds. A simple counterexample is
given by the following  $2\times 2$ matrix $A$ in $\R^2$, $A=\left(
\begin{array}{cc}
                            0 & b \\
                            -a & -1 \\
                          \end{array}
                        \right).$
The eigenvalues of this matrix are $-0.5 \pm i\sqrt{ab-0.25}$, and
if $ab\ge 0.25$, then the spectrum $\sigma(A)$ of $A$ lies in the
half-plane Re$z\leq -0.5$. On the other hand, if, for example, $a=1$
and $b=5$, $u_1=u_2=0.5$, then $(A_Ru,u)>0$.

Inequality Re$(Au,u)\le 0$ means that the operator $A(t)$ is {\it
dissipative}. Such operators often arise in applications (see, e.g.,
\cite{T}). The dissipativity property, defined by the above
inequality,  usually means that the energy in the system is
dissipating, that is, the system is passive. In \cite{R129} a wide
class of passive nonlinear networks is studied, see also
\cite{R118}, Chapter 3.

Our basic results on the stability of the solutions to problem
\eqref{e1}-\eqref{e2} with dissipative operator $A(t)$ are
formulated in Theorems 1.2 and 1.4. Theorem 1.2 contains an auxiliary
result used in the proofs of Theorems 1.2 and 1.4. This result is of
interest by itself and useful in applications.
\begin{theorem}
{\it Assume that Re$(Au,u)\le -|\kappa|\|u\|^2$ for
every $u\in H$ and inequality \eqref{e3} holds. Then the solution to
problem \eqref{e1}$-$\eqref{e2} satisfies an estimate
$\|u(t)\|=O(e^{-(|\kappa|-\epsilon)t})$ as $t\to \infty$. Here
$0<\epsilon< |\kappa|$ can be chosen arbitrarily small if $\|u_0\|$
is sufficiently small.}
\end{theorem}
This theorem implies asymptotic stability in the sense of Lyapunov
of the zero solution to equation \eqref{e1}. Our proof of Theorem 1.2
is new and very short.

We first prove Theorem 1.2 and Theorem 1.4 in Section 2, because the
ideas of our proofs of these theorems are quite similar. Theorem 1.4
contains a new result, and it is not assumed in the formulation of
this theorem that the spectrum of $A(t)$ lies in a half-plane
Re$z\le -|\kappa|$ with $|\kappa|>0$ being a constant independent of
$t$.

Then we prove Theorem 1.3. The result of this theorem is used in the
proofs of Theorems 1.2 and 1.4, and is of general interest. It gives a
bound on solutions to a nonlinear differential inequality. Results
of this type, but considerably less general, were used extensively
in \cite{R499}, where the Dynamical Systems Method (DSM) for solving
operator equations, especially nonlinear equations, was developed.

The ideas of our proofs are quite different from these  in
\cite{DK}.
\begin{theorem}
{\it Let $g(t)\ge 0$ be defined on an interval
$[0,T)$, $T>0$, and have a bounded derivative from the right at
every point of this interval, $\dot{g}(t):=\lim_{s\to
+0}\frac{g(t+s)-g(t)}{s}$. Assume that $g(t)$ satisfies the
following inequality
\begin{equation}
\label{e9}
\dot{g}(t) \le -\gamma(t)g(t) + \alpha(t,g(t)) +
\beta(t),\quad t\in [0,T);\quad g(0)=g_0,
\end{equation}
where $\beta(t)\geq 0$ and $\gamma(t)\geq 0$ are continuous
functions, defined on $[0,\infty)$, and $\alpha(t,v)\geq 0$ is
defined on $[0,\infty)\times [0, \infty)$,
 $\alpha(t,v)$ is non-decreasing as a
function of $v$, locally Lipschitz with respect to $v$,  and
continuous with respect to $t$ on $\R_+:=[0,\infty)$.

 If there exists a function $\mu>0$, continuously differentiable on
 $\R_+$, such that
\begin{equation}
\label{e10}
 \alpha(t,\frac 1{\mu(t)}) + \beta(t)\le \frac 1
{\mu(t)}\left(\gamma(t)-\frac{\dot{\mu}(t)}{\mu(t)}\right),
\qquad\forall
t\ge 0,
\end{equation}
and
\begin{equation}
\label{e11}
g(0)< \frac 1{\mu(0)},
\end{equation}
then $g(t)$ exists for all $t\ge 0$, that is, $g(t)$ can be extended
from $[0,T)$ to $[0,\infty)$,  and  $g(t)$ satisfies  the following
inequality:
\begin{equation}
\label{e12}
0\le g(t)< \frac 1{\mu(t)}, \qquad \forall t\ge 0.
\end{equation}
If $g(0)\le \frac 1{\mu(0)}$, then  $0\le g(t)\le \frac 1{\mu(t)}, 
\qquad \forall t\ge 0$.
}
\end{theorem}
Inequality \eqref{e12} was formulated in \cite{R593} under some
different assumptions, but not proved there. We sketch its  proof
at the end of this paper.

In \cite{R558} inequality \eqref{e9} is studied in the case that
includes $\alpha(t,g)=c_0g^p$, where $p>1$ and $c_0>0$ are
constants, as a particular case: the coefficient $c_0$ in
\cite{R558} was a function of time.

Our second stability result is the following theorem.
\begin{theorem}
{\it Assume that inequality \eqref{e3} holds,
\begin{equation}
\label{e13}
Re(A(t)u,u)\le -\gamma(t) \|u\|^2, \qquad \forall t\ge 0,
\end{equation}
and
\begin{equation}
\label{e14}
\gamma(t)=\frac {c_1}{(1+t)^q}, \quad q\le 1; \quad c_1,q=const>0.
\end{equation}
Suppose that $\epsilon\in (0,c_1)$ is an arbitrary fixed number,
$\lambda= \left(\frac{c_0}{\epsilon}\right)^{1/(p-1)}$,
$(p-1)(c_1-\epsilon)\ge q$,
and  $\|u(0)\|\le \frac {1}{\lambda}$.

Then the unique solution to \eqref{e1}$-$\eqref{e2} exists on all of
$\R_+$ and
\begin{equation}
\label{e15}
0\le \|u(t)\|\le \frac {1}{\lambda(1+t)^{c_1-\epsilon}}.
\end{equation}
}
\end{theorem}

{\bf Remark 1.} {\it One may change the formulation of Theorem 1.4 as
follows: if for some positive constants  $\lambda$ and $\nu>0$
inequalities \eqref{e27} and \eqref{e22} (see below) hold, then
inequality $\|u(t)\|\le \frac 1 {\lambda (1+t)^{\nu}}$ holds for all
$t\ge 0$ for the solution to problem \eqref{e1}$-$\eqref{e2}, as
follows from the proof of Theorem 1.4, given in Section 3.}

The rate of decay of the solution $u(t)$ as $t\to \infty$, obtained
in Theorem 1.4, is not necessarily the best possible.  The result in
Theorem 1.4 is novel and interesting because no assumption of the type
$\gamma(t)\ge \gamma_0>0$, where $\gamma_0$ is a constant, is made.
This allows one to study, for instance, evolution problems with
elliptic operators $A(t)$ the ellipticity constant $\lambda(t)$ of
which may tend to zero as $t\to \infty$. Here $\lambda(t)$ is the
smallest eigenvalue of the matrix $a_{ij}(t)$ of the elliptic
operator $A(t)$. An example is given in Remark 2, at the end of the
paper. 

We have assumed above that $A(t)$ is a bounded linear operator,
since this assumption is basic in the book \cite{DK}, and in the 
Introduction to our paper a comparison was  made with the results in 
\cite{DK}.  However,
boundedness of $A(t)$ was not used in our arguments. If $A(t)$ is a
bounded linear operator satisfying the assumptions of Theorems 1.2 or
1.4, then one can guarantee the global existence of the solution to
evolution problem \eqref{e1}$-$\eqref{e2}. If $A(t)$ is an unbounded
linear operator for which the global existence of $u(t)$ holds,
then our arguments, which lead to estimate \eqref{e15}, remain
valid. In the example given in Remark 2, the operator
$A(t)=\gamma(t)(\Delta -I)$, where $\Delta$ is a selfadjoint
realization of the Laplacian in $H=L^2(R^3)$, and $I$ is the identity 
operator in $H$. For this $A(t)$ one knows that the solution 
$u(t)$ to
problem \eqref{e1}$-$\eqref{e2} exists globally, so Theorem 1.4 is
applicable.

In Section 2 proofs are given.

\section{Proofs}
\label{sec2}

\begin{proof}
{(\it Proof of Theorem 1.2)}.

Multiply \eqref{e1} by $u$, denote $g=g(t):=\|u(t)\|$, take the real
part, and use the assumption \eqref{e13} with
$\gamma(t)=|\kappa|=const>0$, to get
\begin{equation}
\label{e16}
g\dot{g}\le -|\kappa|g^2+c_0g^{p+1}, \qquad p>1.
\end{equation}
If $g(t)>0$ then the derivative $\dot{g}$ does exist, as one can
easily check. If $g(t)=0$ on an open subset of $\R_+$, then the
derivative $\dot{g}$ does exist on this subset and $\dot{g}(t)=0$ on
this subset. If $g(t)=0$ but in any neighborhood $(t-\delta,
t+\delta)$ there are points at which $g$ does not vanish, then by
$\dot{g}$ we understand the derivative from the right, that is,
$$\dot{g}(t):= \lim_{s\to +0}\frac {g(t+s)-g(t)}{s}=\lim_{s\to +0}\frac
{g(t+s)}{s}.$$
This limit does exist and is equal to $\|\dot{u}(t)\|$. Indeed, the
function $u(t)$ is continuously differentiable, so
$$\lim_{s\to +0}\frac {\|u(t+s)\|}{s}=\lim_{s\to +0}
\|\dot{u}(t)+o(1)\|=\|\dot{u}(t)\|.$$
The assumption about the existence of the bounded derivative
$\dot{g}(t)$ from the right in Theorem 1.3 was made because the
function $\|u(t)\|$ does not have, in general, a derivative in the
usual sense at the points $\tau$ at which $\|u(\tau)\|=0$, no matter
how smooth the function $u(t)$ is at the point $\tau$. However, as
we have proved above, the derivative $\dot{g}(t)$ from the right
does exist always, if $u(t)$ is continuously differentiable at the
point $t$.

Since $g\ge 0$, the inequality \eqref{e16}  yields inequality
\eqref{e9} with $\gamma(t)=|\kappa|=const>0$, $\beta(t)=0$, and
$\alpha(t,g)=c_0 g^p$. Inequality \eqref{e10} takes the form
\begin{equation}
\label{e17}
\frac {c_0}{\mu^p(t)}\leq \frac 1 {\mu(t)}\left(|\kappa|-\frac
{\dot{\mu}(t)}{\mu(t)} \right), \qquad  \forall t\ge 0.
\end{equation}
Let
\begin{equation}
\label{e18}
\mu(t)=\lambda e^{bt}, \qquad \lambda,b=const>0,
\end{equation}
and choose the constants $\lambda$ and $b$ later. Then inequality
\eqref{e17} takes the form
\begin{equation}
\label{e19}
\frac {c_0}{\lambda^{p-1}e^{(p-1)bt}} +b\leq  |\kappa|, \qquad \forall
t\ge 0.
\end{equation}
This inequality holds if
\begin{equation}
\label{e20}
\frac {c_0}{\lambda^{p-1}} +b\leq  |\kappa|.
\end{equation}
Let $\epsilon>0$ be an arbitrary small fixed number. Choose
$b=|\kappa|-\epsilon>0$. Then  \eqref{e20} holds if
\begin{equation}
\label{e21}
\lambda\geq \big(\frac{c_0}{\epsilon} \big)^{\frac 1 {p-1}}.
\end{equation}
Condition \eqref{e11} holds if
\begin{equation}
\label{e22}
\|u_0\|=g(0)\le \frac 1 {\lambda}.
\end{equation}
From \eqref{e21}, \eqref{e22} and \eqref{e12} one gets
\begin{equation}
\label{e23}
0\le g(t)=\|u(t)\|\le \frac {e^{-(|\kappa|-\epsilon)t}} {\lambda},\qquad
\forall t\ge 0.
\end{equation}
Theorem 1.2 is proved. 
\end{proof}
\begin{proof}
{(\it Proof of Theorem 1.4.)}

We start with inequality \eqref{e17}, let
\begin{equation}
\label{e24}
\mu(t)=\lambda (1+t)^\nu, \qquad \lambda, \nu=const>0,
\end{equation}
and choose the constants $\lambda$ and $\nu$ later. Inequality
\eqref{e17} holds if
\begin{equation}
\label{e25}
\frac{c_0}{\lambda^{p-1} (1+t)^{(p-1)\nu}} +\frac {\nu}{1+t}\le \frac
{c_1}{(1+t)^q}, \qquad \forall t\ge 0.
\end{equation}
If
\begin{equation}
\label{e26}
q\le 1, \quad (p-1)\nu\ge q,
\end{equation}
then inequality \eqref{e25} holds if
\begin{equation}
\label{e27}
\frac {c_0}{\lambda^{p-1}} +\nu\le c_1.
\end{equation}
Let $\epsilon>0$ be an arbitrary small number. Choose
\begin{equation}
\label{e28}
\nu= c_1-\epsilon.
\end{equation}
Then \eqref{e27}  holds if \eqref{e21} holds. Inequality \eqref{e11}
holds if \eqref{e22} holds. Combining \eqref{e21}, \eqref{e22} and
\eqref{e12}, one obtains
\begin{equation}
\label{e29}
0\le \|u(t)\|=g(t)\le \frac 1 {\lambda (1+t)^{c_1 -\epsilon}}, \qquad
\forall
t\ge 0.
\end{equation}
Choose $\lambda= \big(\frac{c_0}{\epsilon} \big)^{\frac 1 {p-1}}$.
Then inequality \eqref{e27} holds because of \eqref{e28}. Inequality
\eqref{e11} holds because we have assumed in Theorem 1.4 that
$\|u(0)\|\le \frac {1}{\lambda}$. Thus, the desired inequality
\eqref{e15} holds by Theorem 1.3.

Theorem 1.4 is proved. 
\end{proof}
\begin{proof}
{(\it Proof of Theorem 1.3.)}

Define
\begin{equation}
\label{e30}
v(t):=g(t)a(t), \quad a(t):= e^{\int_{t_0}^t\gamma(s)ds}, \quad
\eta(t):=\frac {a(t)}{\mu(t)}.
\end{equation}
Then inequality  \eqref{e9} takes the form
\begin{equation}
\label{e31}
\dot{v}(t) \le a(t)[\alpha\left(t,\frac
{v(t)}{a(t)}\right) + \beta (t)],\qquad v(0)=g(0):= g_0,
\end{equation}
and
\begin{equation}
\label{e32}
\dot{\eta}(t) = \frac{a(t)}{\mu(t)}[\gamma(t)- \frac
{\dot{\mu}(t)}{\mu(t)}].
\end{equation}
From inequalities \eqref{e11} and  \eqref{e10} one gets
\begin{equation}
\label{e33}
v(0)< \frac{1}{\mu(0)}=\eta(0), \quad \dot {v}(0)\leq \dot{\eta}(0). 
\end{equation}
Thus, $v(t)<\eta(t)$ on some interval $[0,T]$.
Inequalities \eqref{e31}, \eqref{e32}, and \eqref{e10} imply
\begin{equation}
\label{e34}
\dot{v}(t) \le \dot{\eta}(t), \quad t\in [0,T].
\end{equation}
It follows from inequalities \eqref{e33} and \eqref{e34}  that
\begin{equation}
\label{e35}
v(t) < \eta(t), \qquad \forall t\ge 0.
\end{equation}
From  inequalities \eqref{e35} and \eqref{e30} one obtains
\begin{equation}
\label{e36}
a(t)g(t)=v(t) < \eta(t)= \frac {a(t)}{\mu(t)}, \qquad  \forall t\ge 0.
\end{equation}
Since $a(t)>0$, inequality \eqref{e36} is equivalent to inequality
\eqref{e12}. This essentially completes the major part of the proof
of inequality \eqref{e12}. The last conclusion of Theorem 1.3 can be
obtained by a standard limiting procedure. 

Let us explain in detail why
inequality \eqref{e36} holds for all $t\ge 0$. The right-hand side
of inequality \eqref{e36} is defined for all  $t\ge 0$. The function
$g(t)$, a solution to inequality \eqref{e9}, exists on every
interval on which $v(t)$ exists, and $v(t)$, the solution to
inequality \eqref{e31}, exists on every interval on which the
solution $w(t)$ to the problem
\begin{equation}
\label{e37}
\dot{w}(t)=a(t)[\alpha(t,\frac {w(t)}{a(t)} + \beta (t)],\qquad
w(0)=v(0)
\end{equation}
exists. It follows from inequality \eqref{e31} and equation
\eqref{e37} that $v(t)\le w(t)$ on every interval $[0,T)$ on which
$w$ exists. We have already proved that the solution to problem
\eqref{e37} (which also is a solution to problem \eqref{e31})
satisfies the estimate
\begin{equation}
\label{e38}
0\le w(t)\le \frac {a(t)}{\mu(t)}
\end{equation}
on every interval on which  $w$ exists. We claim that estimate
\eqref{e38} implies that $w$ exists for all $t\ge 0$, in other
words, that $T=\infty$. Indeed, according to the known result (see,
e.g., \cite{H}, Theorem 3.1 in Chapter 2), if the maximal interval
$[0,T)$ of the existence of the solution to problem \eqref{e37} is
finite, that is $T<\infty$, then $\lim_{t\to T-0}w(t)=\infty$. This,
however, cannot happen because of the inequality  \eqref{e38}, since
the function $\frac {a(t)}{\mu(t)}$ is bounded for every $t\ge 0$.

Theorem 1.3 is proved. 
\end{proof}
{\bf Remark 2.} Let $H=L^2(R^3)$,  $A(t)=\gamma(t) A$, where $A=A^*$
is a selfadjoint operator in $H$ which is the closure of a symmetric
operator $\Delta -I$ with the domain of definition
$C^\infty_0(R^3)$. Here $\Delta$ is the Laplacian. Let $\gamma(t)$
be defined in \eqref{e14} with $c_1=1$, $q=0.5$. let
$\epsilon=0.01$, $p=3$, $c_0=1$, $\lambda=10$, $\nu=0.99$. Assume
that $\|u_0\|\le (0.99)^{-1}$. Theorem 1.4 yields the following
estimate $\|u(t)\|\le 0.1 (1+t)^{-0.99}$ for the solution $u(t)$ to
problem \eqref{e1}$-$\eqref{e2} with the defined above $A(t)$ and a
nonlinearity satisfying condition \eqref{e3}.

\label{lastpage-01}
\end{document}